\documentclass[11pt]{article}

\usepackage[utf8]{inputenc}
\usepackage{mathtools,amsfonts,dsfont,amsthm,mathrsfs,amssymb,stmaryrd}
\usepackage{enumitem}
\usepackage{fullpage}
\usepackage{tikz}
\usepackage{thm-restate}
\usepackage{fullpage}
\usepackage{hyperref}
\usepackage{ifthen}

\usetikzlibrary{fit, backgrounds, matrix, arrows.meta}
\usetikzlibrary{calc,arrows,positioning}
\usetikzlibrary{decorations}
\usetikzlibrary{decorations.pathmorphing}
\usetikzlibrary{decorations.pathreplacing}
\usetikzlibrary{decorations.shapes}
\usetikzlibrary{decorations.text}
\usetikzlibrary{decorations.markings}
\usetikzlibrary{decorations.fractals}
\usetikzlibrary{decorations.footprints}

\definecolor{g-green}{rgb}{0.235, 0.659, 0.322}
\definecolor{g-blue}{rgb}{0.0, 0.5, 1.0}

\newtheorem{theorem}{Theorem}
\newtheorem{lemma}[theorem]{Lemma}
\newtheorem{proposition}[theorem]{Proposition}
\newtheorem{corollary}[theorem]{Corollary}

\theoremstyle{definition}

\newtheorem{problem}[theorem]{Problem}

\newcommand{\dic}{\vec{\chi}}

\newcommand{\bid}{\overleftrightarrow}
\renewcommand{\phi}{\varphi}
\renewcommand{\epsilon}{\varepsilon}

\DeclareMathOperator{\father}{father}
\DeclareMathOperator{\UG}{UG}
\DeclareMathOperator{\dig}{dig}
\DeclareMathOperator{\Mad}{Mad}

\DeclareMathOperator{\tw}{tw}

\newcommand{\ind}[1]{\langle #1 \rangle}

\renewcommand{\epsilon}{\varepsilon}

\renewcommand{\phi}{\varphi}
\let\le\leqslant
\let\ge\geqslant
\let\leq\leqslant
\let\geq\geqslant

\begin{document}

\title{Dichromatic number of chordal graphs}

\author{Stéphane Bessy\footnote{LIRMM, Université de Montpellier, CNRS, Montpellier, France.}
\and
Frédéric Havet\footnote{CNRS, Université Côte d'Azur, I3S, Inria, Sophia-Antipolis, France.}
\and
Lucas Picasarri-Arrieta\footnote{National Institute of Informatics, Tokyo, Japan.\newline Research supported by research grants DIGRAPHS ANR-19-CE48-0013, ANR-17-EURE-0004, and JPMJAP2302.}}

\date{}

\maketitle

\begin{abstract}

The dichromatic number $\dic(D)$ of a digraph $D$ is the minimum integer $k$ such that $D$ admits a $k$-dicolouring, \textit{i.e.} a partition of its vertices into $k$ acyclic subdigraphs. We say that a digraph $D$ is a super-orientation of an undirected graph $G$ if $G$ is the underlying graph of $D$. If $D$ does not contain any pair of symmetric arcs, we just say that $D$ is an orientation of $G$.

In this work, we give both lower and upper bounds on the dichromatic number of super-orientations of chordal graphs. 
In general, the dichromatic number of such digraphs is bounded above by the clique number of the underlying graph (because chordal graphs are perfect). However, this bound can be improved when we restrict the symmetric part of such a digraph.

Let $D=(V,A)$ be a super-orientation of a chordal graph $G$. Let $B(D)$ be the undirected graph with vertex set $V$ in which $uv$ is an edge if and only if both $uv$ and $vu$ belongs to $A$.  An easy greedy procedure shows $\dic(D) \leq \left\lceil \frac{\omega(G) + \Delta(B(D))}{2} \right\rceil$. We show that this bound is best possible by constructing, for every fixed $k,\ell$ with $k\geq \ell+1$, a super-orientation $D_{k,\ell}$ of a chordal graph $G_{k,\ell}$ such that $\omega(G_{k,\ell}) = k$, $\Delta(B(D_{k,\ell})) = \ell$ and $\dic(D_{k,\ell}) = \left\lceil\frac{k+\ell}{2}\right\rceil$.
When $\Delta(B(D)) = 0$ (\textit{i.e.} $D$ is an orientation of $G$), we give another construction showing that this is tight even for orientations of interval graphs.
Next, we show that $\dic(D) \leq \frac{1}{2}\omega(G) + O(\sqrt{d \cdot \omega(G)})$ with $d$ the maximum average degree of $B(D)$.
Finally, we show that if $B(D)$ contains no $C_4$ as a subgraph, then $\dic(D) \leq \left\lceil \frac{\omega(G)+3}{2} \right\rceil$. We justify that this is almost best possible by constructing, for every fixed $k$, a super-orientation $D_{k}$ of a chordal graph $G_{k}$ with clique number $k$ such that $B(D_k)$ is a disjoint union of paths and $\dic(D_k) = \left\lfloor \frac{k+3}{2} \right\rfloor$.

We also exhibit a family of orientations of cographs for which the dichromatic number is equal to the clique number of the underlying graph.
\end{abstract}

\section{Introduction}
We denote by $[k]$ the set $\{1, \dots , k\}$. 
 Given an undirected graph $G=(V,E)$ and a positive integer $k$, a \textit{$k$-colouring} of $G$ is a function $\alpha : V \xrightarrow{} [k]$. It is \textit{proper} if, for every edge $xy\in E$, we have $\alpha(x) \neq \alpha(y)$. So, for every $i\in [k]$, $\alpha^{-1}(i)$ induces an independent set on $G$. The \textit{chromatic number} of $G$, denoted by $\chi(G)$, is the smallest $k$ such that $G$ admits a proper $k$-colouring. An undirected graph is {\it chordal} if it does not contain any induced cycle of length at least 4. Proper colourings of chordal graphs have been largely studied and it is well-known that chordal graphs are perfect. Recall that a graph $G$ is {\it perfect} if every induced subgraph $H$ of $G$ satisfies $\chi(H) = \omega(H)$, where $\omega(H)$ denotes the size of a largest clique in $H$.
 
\medskip

We refer the reader to~\cite{bang2009} for notation and terminology on digraphs not explicitly defined in this paper. Let $D=(V,A)$ be a digraph. A \textit{digon} is a pair of arcs in opposite directions between the same vertices. A \textit{simple arc} is an arc which is not in a digon. 
An \textit{oriented graph} is a digraph with no digon. The \textit{bidirected graph} associated with a graph $G$, denoted by $\bid{G}$,  is the digraph obtained from $G$ by replacing every edge by a digon. The \textit{underlying graph} of $D$, denoted by $\UG(D)$, is the undirected graph with vertex set $V(D)$ in which $uv$ is an edge if and only if $uv$ or $vu$ is an arc of $D$. We say that $D$ is a \textit{super-orientation} of $\UG(D)$, and it is an \textit{orientation} of $\UG(D)$ if $D$ is an oriented graph. A \textit{tournament} on $n$ vertices is an orientation of the complete graph on $n$ vertices. The \textit{bidirected graph} of $D$, denoted by $B(D)$, is the undirected graph $G$ with vertex set $V(D)$ in which $uv$ is an edge if and only if $uv$ is a digon of $D$. We denote by $\bid{\omega}(D)$ the size of a largest bidirected clique of $D$, \textit{i.e.} the size of the largest clique of $B(D)$.

\medskip

In 1982, Neumann-Lara~\cite{neumannlaraJCT33} introduced the notions of dicolouring and dichromatic number, which generalize the ones of proper colouring and chromatic number. 
For a positive integer $k$, a \textit{$k$-colouring} of $D=(V,A)$ is a function $\alpha : V \xrightarrow{} [k]$. It is a \textit{$k$-dicolouring} if $\alpha^{-1}(i)$ induces an acyclic subdigraph in $D$ for each $i \in [k]$. In other words, no directed cycle of $D$ is monochromatic in $\alpha$. The \textit{dichromatic number} of $D$, denoted by $\dic(D)$, is the smallest $k$ such that $D$ admits a $k$-dicolouring. 

There is a one-to-one correspondence between the proper $k$-colourings of a graph $G$ and the $k$-dicolourings of its associated bidirected graph $\bid{G}$, and in particular $\chi(G) = \dic(\bid{G})$. Hence every result on proper colouring of undirected graphs can be seen as a result on dicolouring of bidirected graphs, and it is natural to study whether the result can be extended to all digraphs. 
Indeed, a lot of classical results on graph proper colourings have already been extended to digraphs dicolouring. For instance, Brooks' Theorem (Brooks~\cite{brooksMPCPS37}) has been generalised to digraphs by Harutyunyan and Mohar in~\cite{harutyunyanJDM25} (see also~\cite{aboulkerDM113193}). Another example is the celebrated Strong Perfect Graph Theorem (Chudnovsky, Robertson, Seymour and Thomas~\cite{chudnovskyAM164}) extended to digraphs by Andres and Hochst{\"a}ttler in~\cite{andresJGT79} (the proof is strongly based on the result of Chudnovsky et al). A digraph $D$ is \textit{perfect} if $\dic(H) = \bid{\omega}(H)$ for every induced subdigraph $H$ of $D$.

\begin{theorem}[Andres and Hochst{\"a}ttler~\cite{andresJGT79}]
    \label{thm:perfect_digraphs}
    A digraph $D$ is perfect if and only if $B(D)$ is perfect and $D$ does not contain an induced directed cycle of length at least 3.
\end{theorem}

We refer the interested reader to~\cite{meisterTCS463}, in which the authors define a class of chordal digraphs, which extends the class of undirected chordal graphs. One can easily prove that every digraph $D$ in this class is actually a perfect digraph, so it satisfies $\dic(D) = \bid{\omega}(D)$ by Theorem~\ref{thm:perfect_digraphs}.

\medskip

In this work, we look for lower and upper bounds on the dichromatic number of orientations and super-orientations of chordal graphs. Dicolourings of such digraphs have also been studied in~\cite{aboulkerSIAM36}, in which the authors characterise exactly the digraphs $H$ for which there exists $c_H\in \mathbb{N}$ such that every oriented chordal graph $\Vec{G}$ with $\dic(\vec{G}) \geq c_H +1$ contains $H$ as an induced subdigraph.

The very first interesting class of such digraphs are tournaments for which the question has been settled by Erd{\"o}s, Gimbel and Kratsch in~\cite{erdosJGT15}. They showed that the dichromatic number of a tournament $T$ on $n$ vertices is always at most $O\left(\frac{n}{\log n}\right)$, and that this bound is tight (up to a constant factor). One can ask if this result is true not only for tournaments but for all orientations of chordal graphs. That is, do we always have $\dic(\Vec{G}) = O\left(\frac{\omega(G)}{\log \omega(G)}\right)$ when $\Vec{G}$ is an orientation of a chordal graph $G$~? We answer this by the negative. Indeed, we show in Section~\ref{section:orientations_interval} that it is not even true for orientations of interval graphs. Recall that an {\it interval graph} is obtained from a set of intervals on the real line: the intervals are the vertices and there is an edge between two intervals if and only if they intersect. It is well-known that interval graphs are chordal.

\begin{restatable}{theorem}{constructioninterval}
    \label{thm:construction_interval}
    For every fixed $k \in \mathbb{N}$, there exists an interval graph $G_k$ and an orientation $\vec{G}_k$ of this graph such that $\omega(G_k) = k$ and $\dic(\Vec{G}_k) \geq \lceil \frac{k}{2} \rceil$.
\end{restatable}

On the positive side, if $\Vec{G}$ is the orientation of a proper interval graph $G$ (which is an interval graph where each interval has length exactly one), then $\dic(\Vec{G}) = O\left(\frac{\omega(G)}{\log(\omega(G))}\right)$, as proved in~\cite{aboulkerSIAM36}. The key idea is that $G$ admits a partition $(V_1,V_2)$ of its vertex-set such that both $G\ind{V_1}$ and $G\ind{V_2}$ are disjoint union of cliques.

\medskip 

Another well-known class of perfect graphs is the one of cographs. The \textit{join} of two undirected graphs $G_1$ and $G_2$ is the graph built from the disjoint union of $G_1$ and $G_2$ where every edge between vertices of $G_1$ and vertices of $G_2$ are added. Cographs form the smallest class of graphs containing the single-vertex graph that is closed under disjoint union and the join operation. One can easily prove that the oriented graphs built in the proof of Theorem~\ref{thm:construction_interval} are indeed orientations of cographs. In Section~\ref{section:orientations_cograph}, we improve this result for cographs in general.

\begin{restatable}{theorem}{constructioncograph}
    \label{thm:construction_cograph}
    For every $k \in \mathbb{N}$, there exists an orientation  $\vec{G}_k$ of a cograph $G_k$ such that $\dic(\Vec{G}_k) = \omega(G_k) = k$.
\end{restatable}

\medskip

Next we consider super-orientations of chordal graphs. If $D$ is a super-orientation of a chordal graph $G$, then obviously $\dic(D) \leq \omega(G)$ because $\dic(D) \leq \chi(G) = \omega(G)$. Note that we cannot expect any improvement of this bound in general, because if $D$ is the bidirected graph $\bid{G}$ then $\dic(D) = \omega(G)$. But one can ask what happens if we restrict the structure of $B(D)$, the bidirected graph of $D$. 

In Section~\ref{section:bounded_degree}, we consider digraphs for which the bidirected graph has bounded maximum degree. Using the degeneracy of the underlying graph, we show the following easy proposition.
\begin{restatable}{proposition}{propboundedmaxdegree}
    \label{prop:bounded_max_degree}
    Let $D$ be a super-orientation of a chordal graph $G$. Then 
    \[\dic(D) \leq \left\lceil \frac{\omega(G) + \Delta(B(D))}{2} \right\rceil .\] 
\end{restatable}

This proposition is best possible when $\Delta(B(D)) = 0$ by Theorem~\ref{thm:construction_interval}. In the following, we show that it is indeed best possible for every fixed value of $\Delta(B(D))$.
\begin{restatable}{theorem}{constructionchordal}
    \label{construction_chordal}
    For every fixed $k,\ell \in \mathbb{N}$ such that $k\geq \ell+1$, there exists a chordal graph $G_{k,\ell}$ and a super-orientation $D_{k,\ell}$ of $G_{k,\ell}$ such that $\omega(G_{k,\ell}) = k$, $\Delta(B(D_{k,\ell})) = \ell$ and $\dic(D_{k,\ell}) = \left\lceil \frac{k + \ell}{2} \right\rceil$.
\end{restatable}

The \textit{maximum average degree} of an undirected graph $G$ is $\Mad(G) = \max \left\{ \frac{2|E(H)|}{|V(H)|} \mid H \text{ subgraph of } G\right\}$.
In Section~\ref{section:bounded_mad}, we show the following bound on digraphs $D$ for which $\Mad(B(D))$ is bounded.
\begin{restatable}{theorem}{boundedmad}
    \label{thm:bounded_mad}
    Let $D$ be a super-orientation of a chordal graph $G$. If $\Mad(B(D)) \leq d$, then
    \[
        \dic(D) \leq \frac{1}{2}\omega(G) + O(\sqrt{d\cdot \omega(G)}). 
    \]
\end{restatable}

Finally in Section~\ref{section:no_C4} we show the following bound on super-orientations $D$ of chordal graphs that do not contain $\bid{C_4}$. We say that a graph $G$ is $C_4$-free if $G$ does not contain $C_4$ as a subgraph.
\begin{restatable}{theorem}{Cfourfree}
    \label{thm:C4_free}
    Let $D$ be a super-orientation of a chordal graph $G$. If $B(D)$ is $C_4$-free, then 
    \[\dic(D) \leq \left\lceil \frac{\omega(G)+3}{2} \right\rceil.\]
\end{restatable}

We also prove that the bound of Theorem~\ref{thm:C4_free} is almost tight by proving the following.
\begin{restatable}{theorem}{constructionCfourfree}
    \label{construction:C4_free}
    For every fixed $k\geq 3$ and every $n\geq \mathbb{N}$, there exists a super-orientation $D_{k,n}$ of a chordal graph $G_{k,n}$ on at least $n$ vertices such that $B(D_{k,n})$ is a disjoint union of paths, $\omega(G_{k,n}) = k$ and $\dic(D_{k,n}) = \left\lfloor \frac{k+3}{2} \right\rfloor$.
\end{restatable}

A \textit{tree-decomposition} of a graph $G=(V,E)$ is a pair $(T,\mathcal{X})$ where $T=(I,F)$ is a tree, and $\mathcal{X}=(B_i)_{i\in I}$ is a family of subsets of $V(G)$, called \textit{bags} and indexed by the vertices of $T$, such that:
\begin{enumerate}
    \item each vertex $v\in V$ appears in at least one bag, \textit{i.e.} $\bigcup_{i\in I} B_i= V$,
    \item for each edge $e = xy \in E$, there is an $i\in I$ such that $x,y \in B_i$, and 
    \item for each $v\in V$, the set of nodes indexed by $\{ i \mid i\in I, v\in B_i\}$ forms a subtree of $T$.
\end{enumerate}
The \textit{width} of a tree decomposition is defined as $\max_{i\in I} \{|B_i| -1\}$. The \textit{treewidth} of $G$, denoted by $\tw (G)$, is the minimum width of a tree-decomposition of $G$.
It is well-know that every graph $G$ is a subgraph of a chordal graph $G'$ with $\omega(G') = \tw(G) + 1$. Hence the following is a direct consequence of Proposition~\ref{prop:bounded_max_degree} and Theorems~\ref{thm:bounded_mad} and~\ref{thm:C4_free}.
\begin{corollary}
    \label{cor-tw}
    Let $D$ be a super-orientation of $G$. Then we have:
    \begin{itemize}
        \item $\dic(D) \leq \left\lceil \frac{\tw(G) + \Delta(B(D)) + 1}{2} \right\rceil$, and
        \item $\dic(D) \leq \frac{1}{2}\tw(G) + O(\sqrt{\Mad(B(D)) \cdot \tw(G)})$, and
        \item $\dic(D) \leq \left\lceil \frac{\tw(G)+4}{2} \right\rceil$ if $B(D)$ is $C_4$-free.
    \end{itemize}
\end{corollary}

\section{Definitions and preliminary results}

Let $G=(V,E)$ be an undirected graph. A \textit{perfect elimination ordering} of $G$ is an ordering $v_1,\dots,v_n$ of its vertex-set such that, for every $i\in [n]$, the subgraph of $G$ induced by $N(v_i) \cap \{v_{i+1},\dots,v_n\}$ is a clique.
\begin{proposition}[Folklore]
    \label{prop:perfect_ordering}
    A graph $G$ is chordal if and only if $G$ admits a perfect elimination ordering.
\end{proposition}

\begin{proposition}[Folklore]
    \label{prop:treewidth_chordal}
    The treewidth of a chordal graph $G$ is exactly $\omega(G) - 1$.
\end{proposition}

A tree-decomposition $(T,{\cal X})$ is {\it reduced} if, for every $tt' \in E(T)$, $X_t \setminus X_{t'}$ and $X_{t'} \setminus X_t$ are non-empty. It is easy to see that any graph $G$ admits an optimal (i.e., of width $\tw(G)$) tree-decomposition which is reduced (indeed, if $X_t \subseteq X_{t'}$ for some edge  $tt' \in E(T)$, then contract this edge and remove $X_t$ from ${\cal X}$). 

A tree-decomposition $(T,{\cal X})$ of a graph $G$ of width $k\geq 0$ is {\it full} if every bag has size exactly $k+1$. It is {\it valid} if $|X_t \setminus X_{t'}|=|X_{t'}\setminus X_t|=1$ for every $tt' \in E(T)$. Note that any valid tree-decomposition is full and reduced. 

The following result is well-known, see for instance~\cite{bodlaenderSIAMJC25}. We give here a short proof for sake of completeness.

\begin{lemma}\label{lem:valid_decomposition}
Every graph $G=(V,E)$ admits a valid tree-decomposition of width $\tw(G)$.
\end{lemma}
\begin{proof}
Let $(T,{\cal X})$ be an optimal reduced tree-decomposition of $G=(V,E)$, which exists by the remark above the lemma. We will progressively modify $(T,{\cal X})$ in order to make it first full and then valid.

While the current decomposition is not full, let $tt' \in E(T)$ such that $|X_t| < |X_{t'}|=\tw(G)+1$ and let $v \in X_{t'} \setminus X_t$. Add $v$ to $X_t$. The obtained decomposition is still a tree-decomposition. Moreover, the updated decomposition remains reduced all along the process, as since $|X_t| < |X_{t'}|$ and the initial decomposition is reduced, $X_{t'}$ must contain another vertex $u \neq v$ with $u \notin X_t$. At the end of the process, we obtain an optimal decomposition $(T,{\cal X})$ that is full. 

Now, while $(T,{\cal X})$ is not valid, let $tt' \in E(T)$, $x,y \in X_t \setminus X_{t'}$ and $u,v \in X_{t'} \setminus X_t$ (such an edge of $T$ and four distinct vertices of $V$ must exist since  $(T,{\cal X})$ is full and reduced but not valid). Then, add a new node $t''$ to $T$, with corresponding bag $X_{t''}=(X_{t'} \setminus \{u\}) \cup \{x\}$ and replace the edge $tt'$ in $T$ by the two edges $tt''$ and $t''t'$. Clearly, subdividing the edge $tt'$ by adding a bag $X_{t''}=X_{t'}\setminus\{u\} \cup \{x\}$ still leads to an optimal full tree-decomposition of the same width.

Note that, after the application of each step as described above, either the maximum of $|X_t\setminus X_{t'}|$ over all edges $tt' \in E(T)$, or the number of edges $tt' \in E(T)$ that maximize $|X_t\setminus X_{t'}|$, strictly decreases, and none of these two quantities increases. Therefore, the process terminates, and eventually $(T,{\cal X})$ becomes an optimal valid tree-decomposition.
\end{proof}

\medskip

Let $D_1$ and $D_2$ be two digraphs. Let $u_1v_1$ be an arc of $D_1$ and $v_2u_2$ be an arc of $D_2$. The \textit{directed Haj\'os join} of $D_1$ and $D_2$, denoted by $D_1\triangledown D_2$, is the digraph obtained from the union $D_1\cup D_2$ by deleting the arcs $u_1v_1$ as well as $v_2u_2$, identifying the vertices $v_1$ and $v_2$ into a new vertex $v$ and adding the arc $u_1u_2$.
\begin{theorem}[Bang-Jensen et al.~\cite{bangjensenEJC27} (see also~\cite{hoshinoCombinatorica35})]
    \label{thm:hajos_join}
    Let $D_1$ and $D_2$ be two digraphs, then
    \[ \dic(D_1 \triangledown D_2) \geq \min \{ \dic(D_1),\dic(D_2)\}. \]
\end{theorem}

\section{Orientations of interval graphs with large dichromatic number}
\label{section:orientations_interval}

This section is devoted to the proof of Theorem~\ref{thm:construction_interval}.

\constructioninterval*

\begin{proof}
    Let us fix $k\in \mathbb{N}$, we will build an orientation $D_k$ of an interval graph $G_k$ such that $\omega(G_k) = k$ and $\dic(D_k) \geq \left \lfloor \frac{k+1}{2} \right \rfloor$.

We start from one interval $I_1^1$. Then, for every $i$ from $2$ to $k$, we do the following: for each interval $I_{i-1}^s$  we added at step $i-1$, we add $2^{i-1}$ new pairwise disjoint intervals whose union is included in $I_{i-1}^s$, and we associate to each of these new intervals $I_i^\ell$ a distinct binary number $b_{i}^\ell$ on $i-1$ bits. By construction, every new interval intersects exactly $i-1$ other intervals (one for each step). 

Let $G_k$ be the interval graph made of the intervals built above. By construction, $\omega(G_k) = k$. Now we consider $D_k$ the orientation of $G_k$ defined as follows. For every pair $j<i$, we orient the edge $I_j^sI_i^\ell$  from  $I_i^\ell$ to $I_j^s$ if the $j^{\text{th}}$ bit of $b_i^\ell$ is 1, and from $I_j^s$ to $I_i^\ell$ otherwise. Figure~\ref{fig:interval-bound} illustrates the construction of $D_3$. 

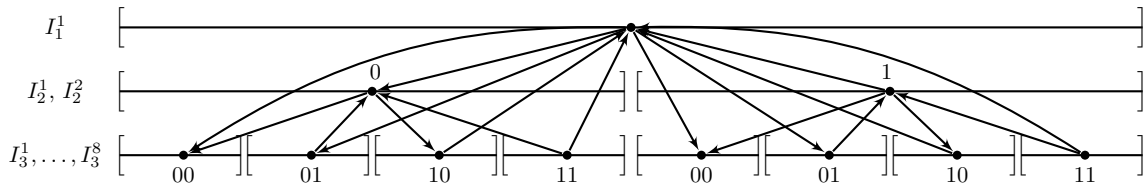
\begin{figure}[hbtp]
  \begin{minipage}{\linewidth}
    \begin{center}	
      \begin{tikzpicture}[thick,scale=.85, every node/.style={transform shape}]
      
	  	\tikzset{edge/.style = {->,> = latex'}}
	    \tikzset{vertex/.style = {circle,fill=black,minimum size=4pt,
                            inner sep=0pt}}
                            
        \draw[] (0,0) -- (16,0);
        \node[] at (16,0) {$\Big]$};
        \node[] at (0,0) {$\Big[$};
        \node[vertex] (0) at (8,0) {};
        
        \draw[] (0,-1) -- (7.9,-1);
        \node[] at (7.9,-1) {$\Big]$};
        \node[] at (0,-1) {$\Big[$};
        \node[] at (4,-0.7) {$0$};
        \node[vertex] (1) at (3.95,-1) {};
        
        \draw[] (8.1,-1) -- (16,-1);
        \node[] at (16,-1) {$\Big]$};
        \node[] at (8.1,-1) {$\Big[$};
        \node[] at (12,-0.7) {$1$};
        \node[vertex] (2) at (12.05,-1) {};
        
        \foreach \i in {0,...,3}{
            \draw[] (0 + \i * 2,-2) -- (1.9 + \i * 2,-2);
            \node[] at (1.9 + \i * 2,-2) {$\Big]$};
            \node[] at (0 + \i * 2,-2) {$\Big[$};
            \node[vertex] (3_\i) at (1 + \i * 2,-2) {};
            
            \draw[] (8.1 + \i * 2,-2) -- (10 + \i * 2,-2);
            \node[] at (10 + \i * 2,-2) {$\Big]$};
            \node[] at (8.1 + \i * 2,-2) {$\Big[$};
            \node[vertex] (4_\i) at (9.1 + \i * 2,-2) {};
        }
        
        \node[] at (1,-2.3) {$00$};
        \node[] at (3,-2.3) {$01$};
        \node[] at (5,-2.3) {$10$};
        \node[] at (7,-2.3) {$11$};
        \node[] at (9.1,-2.3) {$00$};
        \node[] at (11.1,-2.3) {$01$};
        \node[] at (13.1,-2.3) {$10$};
        \node[] at (15.1,-2.3) {$11$};
        
        \draw[edge] (0) to (1);
        \draw[edge] (2) to (0);
        
        \draw[edge] (1) to (3_0);
        \draw[edge] (1) to (3_2);
        \draw[edge] (3_1) to (1);
        \draw[edge] (3_3) to (1);
        
        \draw[edge] (2) to (4_0);
        \draw[edge] (2) to (4_2);
        \draw[edge] (4_1) to (2);
        \draw[edge] (4_3) to (2);
        
        \draw[edge, bend right=17] (0) to (3_0);
        \draw[edge] (0) to (3_1);
        \draw[edge] (3_2) to (0);
        \draw[edge] (3_3) to (0);
        
        \draw[edge] (0) to (4_0);
        \draw[edge] (0) to (4_1);
        \draw[edge] (4_2) to (0);
        \draw[edge, bend right=17] (4_3) to (0);

        \draw (-1,0) node[]{$I_1^1$};
        \draw (-1,-1) node[]{$I_2^1$, $I_2^2$};
        \draw (-1,-2) node[]{$I_3^1,\dots , I_3^8$};
        
      \end{tikzpicture}
      
      \caption{The oriented interval graph $D_3$ (bits of $b_i^\ell$ are read from left to right).}
      \label{fig:interval-bound}
    \end{center}
  \end{minipage}
\end{figure}

Let us prove that $\dic(D_k) \geq \lceil \frac{k}{2} \rceil$. To do this, let $\phi$ be any optimal dicolouring of $D_k$. We will find a tournament $T$ of size $k$ in $D_k$ such that, for each colour $c$ in $\phi$, $c$ appears at most twice in $T$. This will prove that $\phi$ uses at least $\lceil \frac{k}{2} \rceil$ colours, implying the result.

Start from the universal vertex $I_1^1$. Then, for $i\in \{2,\dots, k\}$, we do the following : let $I_{i-1}^s$ be the last vertex added to $T$, we will extend $T$ with a vertex $I_{i}^\ell$ so that $I_{i}^\ell \subseteq I_{i-1}^s$. For each colour $c \in \phi$ that appears exactly twice in $T$, let $x_cy_c$ be a monochromatic arc of $T$ coloured $c$. Then we choose $I_{i}^\ell$ so for each such colour $c$, $x_cy_cI_{i}^\ell$ is a directed triangle. The existence of $I_{i}^\ell$ is guaranteed by construction. This implies that the colour of $I_{i}^\ell$ in $\phi$ appears at most twice in $T$.
\end{proof}

\section{Orientations of cographs with large dichromatic number}
\label{section:orientations_cograph}

This section is devoted to the proof of Theorem~\ref{thm:construction_cograph}.

\constructioncograph*

\begin{proof}
    We define $\vec{G}_1$ as the only orientation of $G_1$, the graph on one vertex. We obviously have $\dic(\vec{G}_1) = \omega(G_1) = 1$, and $G_1$ is a cograph.

    Let us fix $k\geq 1$, we build $\vec{G}_{k+1}$ from $\vec{G}_k$ as follows. Start from $k+1$ disjoint copies $\vec{G}_k^1,\dots,\vec{G}_k^{k+1}$ of $\vec{G}_k$ and $k+1$ new vertices $v_1,\dots,v_{k+1}$. Then, for every $i\in [k+1]$, we add all arcs from $v_i$ to $V(\vec{G}_k^i)$ and all arc from $\bigcup_{j\neq i}V(\vec{G}_k^j)$ to $v_i$. 
    Let $\vec{G}_{k+1}$ be the obtained oriented graph and $G_{k+1}$ be its underlying graph. Figure~\ref{fig:construction_cograph} illustrates the construction of $\vec{G}_3$.

    \begin{figure}[hbtp]
  \begin{minipage}{\linewidth}
    \begin{center}	
      \begin{tikzpicture}[thick,scale=.85, every node/.style={transform shape}]
        \tikzset{edge/.style = {->,> = latex}}
        \tikzset{vertex/.style = {circle,fill=black,minimum size=6pt,
                        inner sep=0pt}}
        \tikzset{bigvertex/.style = {shape=circle,draw,dotted, minimum size=11ex}}
        
        \node[vertex,orange] (d11) at  (-1,0) {};
        \node[] (d1) at  (-1,-1) {$\Vec{G}_1$};
        
        \draw[dotted, line width=1pt] (1,0) circle (1.6ex);
        \draw[dotted, line width=1pt] (1,1) circle (1.6ex);
        \node[vertex,orange] (d3_1) at  (1,0) {};
        \node[vertex,orange] (d3_2) at  (1,1) {};
        \node[vertex,g-blue] (d3_3) at  (2,0) {};
        \node[vertex,g-blue] (d3_4) at  (2,1) {};
        \node[] (d3) at  (1.5,-1) {$\Vec{G}_2$};
        \draw[edge] (d3_3) to (d3_1);
        \draw[edge] (d3_1) to (d3_4);
        \draw[edge] (d3_2) to (d3_3);
        \draw[edge] (d3_4) to (d3_2);
        
        \foreach \i in {0,...,2}{
            \node[vertex,orange] (d41\i) at  (4,0+2*\i) {};
            \node[vertex,orange] (d42\i) at  (4,1+2*\i) {};
            \node[vertex,g-blue] (d43\i) at  (5,0+2*\i) {};
            \node[vertex,g-blue] (d44\i) at  (5,1+2*\i) {};
            \node[bigvertex] (d46\i) at (4.5,0.5+2*\i) {};
            \draw[edge] (d43\i) to (d41\i);
            \draw[edge] (d41\i) to (d44\i);
            \draw[edge] (d42\i) to (d43\i);
            \draw[edge] (d44\i) to (d42\i);
            
            \node[vertex,g-green] (d45\i) at  (8,0.5+2*\i) {};
        }
        \foreach \i in {0,...,2}{
            \foreach \j in {0,...,2}{
                \ifthenelse{\i = \j}
                    {
                        \draw[edge] (d45\i) to (d46\i);
                    }
                    {
                        \draw[edge] (d46\j) to (d45\i);
                    }
                  ;
            }
        }

        \node[] (d3) at  (6,-1) {$\Vec{G}_3$};
        
      \end{tikzpicture}
      \caption{The oriented graphs $\Vec{G}_1$, $\Vec{G}_2$ and $\Vec{G}_3$.}
      \label{fig:construction_cograph}
    \end{center}
  \end{minipage}
\end{figure}
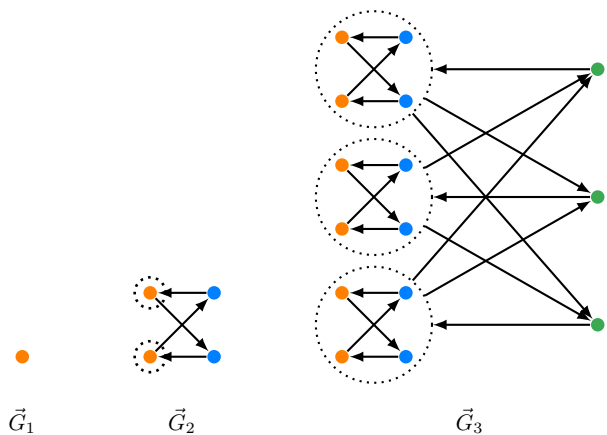
    
    Note first that $G_{k+1}$ is a cograph: the disjoint union of $G_k^1,\dots,G_k^k$ is a cograph, the independent set $v_1,\dots,v_k$ is a cograph, and $G_{k+1}$ is the join of these two cographs. 
    Let us prove by induction on $k$ that $\dic(\Vec{G}_{k}) = \omega(G_{k}) = k$. For $k=1$, the result is immediate, and assume it holds for $k\ge 1$.
    Note first that $\omega(G_{k+1}) = k+1$ since every clique of $G_{k+1}$ contains at most one vertex of $\{v_1,\dots,v_{k+1}\}$ and do not contain two vertices from distinct copies of $G_{k}$. So every maximum clique of $G_{k+1}$ is made of a maximum clique of $G_k$ and one additional vertex $v_i$.

Moreover $\dic(\Vec{G}_{k+1}) \leq \chi(G_{k+1}) = \omega(G_{k+1}) = k+1$.
    Let us now show that the dichromatic number of $\vec{G}_{k+1}$ is at least $k+1$. Assume for the purpose of contradiction that $\vec{G}_{k+1}$ admits a $k$-dicolouring $\phi$. Then there exist $i\neq j$ such that $\phi(v_i) = \phi(v_j)$. Since $\dic(\Vec{G}_k) \geq k$, there exist $x\in V(\vec{G}_k^i)$ and $y\in V(\vec{G}_k^j)$ such that $\phi(x) = \phi(y) = \phi(v_i)=\phi(v_j)$. Hence $v_ixv_jyv_i$ is a monochromatic $\Vec{C}_4$ of $\vec{G}_{k+1}$ coloured with $\phi$, a contradiction.
\end{proof}

\section{Super-orientations of chordal graphs with a bidirected graph having bounded maximum degree}
\label{section:bounded_degree}

This section is devoted to the proofs of Proposition~\ref{prop:bounded_max_degree} and Theorem~\ref{construction_chordal}.

\propboundedmaxdegree*

\begin{proof}
    Let $v_1,\dots,v_n$ be a perfect elimination ordering of $G$ (which exists by Proposition~\ref{prop:perfect_ordering}). Then, in $G$, every vertex $v_i$ has at most $\omega(G)-1$ neighbours in $\{v_{i+1},\dots,v_{n}\}$. Hence, in $D\ind{\{v_i,\dots, v_n\}}$, $d^+(v_i) + d^-(v_i) \leq \omega(G) - 1 + \Delta(B(D))$.

    Thus, considering the vertices from $v_n$ to $v_1$, we can greedily find a dicolouring of $D$ using at most $\left\lceil \frac{\omega(G)+\Delta(B(D))}{2} \right\rceil$ by choosing for $v_i$ a colour that is not appearing in $N^+(v_i) \cap \{v_{i+1},\dots,v_n\}$ or in $N^-(v_i) \cap \{v_{i+1},\dots,v_n\}$.
\end{proof}

\constructionchordal*

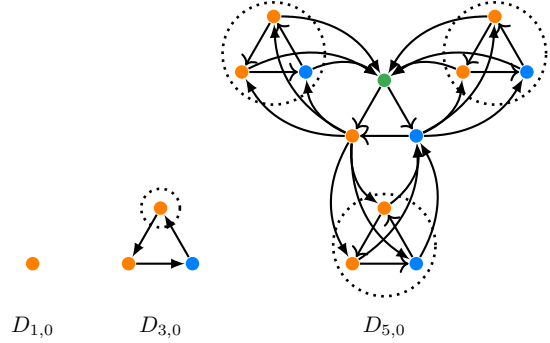
\begin{figure}[hbtp]
  \begin{minipage}{\linewidth}
    \begin{center}	
      \begin{tikzpicture}[thick,scale=.85, every node/.style={transform shape}]
        
  	\tikzset{edge/.style = {->,> = latex}}
        \tikzset{vertex/.style = {circle,fill=black,minimum size=6pt,
                            inner sep=0pt}}
        
        \node[vertex,orange] (d11) at  (-0.50,0) {};
        \node[] (d1) at  (-0.5,-1) {$D_{1,0}$};
        
        \draw[dotted, line width=1pt] (1.5,0.87) circle (2ex);
        \node[vertex,orange] (d3_1) at  (1,0) {};
        \node[vertex,g-blue] (d3_2) at  (2,0) {};
        \node[vertex,orange] (d3_3) at  (1.5,0.87) {};
        \node[] (d3) at  (1.5,-1) {$D_{3,0}$};
        \draw[edge] (d3_1) to (d3_2);
        \draw[edge] (d3_2) to (d3_3);
        \draw[edge] (d3_3) to (d3_1);
        
        \node[vertex,orange] (d5_1) at  (4.5,2) {};
        \node[vertex,g-blue] (d5_2) at  (5.5,2) {};
        \node[vertex,g-green] (d5_3) at  (5,2.87) {};
        \node[] (d5) at  (5,-1) {$D_{5,0}$};
        \draw[->] (d5_2) -- (d5_1) node [midway, auto] () {};
        \draw[->] (d5_3) -- (d5_2) node [midway, auto] () {};
        \draw[->] (d5_3) -- (d5_1) node [midway, auto] () {};
        
        \draw[dotted, line width=1pt] (5-1.73,3.29) circle (0.8);
        \node[vertex,orange] (d5_4) at  (4.5-1.73,3) {};
        \node[vertex,g-blue] (d5_5) at  (5.5-1.73,3) {};
        \node[vertex,orange] (d5_6) at  (5-1.73,3.87) {};
        \draw[->] (d5_4) -- (d5_5) node [midway, auto] () {};
        \draw[->] (d5_5) -- (d5_6) node [midway, auto] () {};
        \draw[->] (d5_6) -- (d5_4) node [midway, auto] () {};
        
        \draw[dotted, line width=1pt] (5+1.73,3.29) circle (0.8);
        \node[vertex,orange] (d5_7) at  (4.5+1.73,3) {};
        \node[vertex,g-blue] (d5_8) at  (5.5+1.73,3) {};
        \node[vertex,orange] (d5_9) at  (5+1.73,3.87) {};
        \draw[->] (d5_7) -- (d5_8) node [midway, auto] () {};
        \draw[->] (d5_8) -- (d5_9) node [midway, auto] () {};
        \draw[->] (d5_9) -- (d5_7) node [midway, auto] () {};
        
        \draw[dotted, line width=1pt] (5,0.29) circle (0.8);
        \node[vertex,orange] (d5_10) at  (4.5,0) {};
        \node[vertex,g-blue] (d5_11) at  (5.5,0) {};
        \node[vertex,orange] (d5_12) at  (5,0.87) {};
        \draw[->] (d5_10) -- (d5_11) node [midway, auto] () {};
        \draw[->] (d5_11) -- (d5_12) node [midway, auto] () {};
        \draw[->] (d5_12)-- (d5_10) node [midway, auto] () {};
        
        \draw[bend left=30, edge] (d5_1) to (d5_4);
        \draw[bend left=30, edge] (d5_1) to (d5_5);
        \draw[bend left=30, edge] (d5_1) to (d5_6);
        \draw[bend left=30, edge] (d5_4) to (d5_3);
        \draw[bend left=30, edge] (d5_5) to (d5_3);
        \draw[bend left=30, edge] (d5_6) to (d5_3);
        \draw[bend right=30, edge] (d5_1) to (d5_10);
        \draw[bend right=30, edge] (d5_1) to (d5_11);
        \draw[bend right=30, edge] (d5_1) to (d5_12);
        \draw[bend right=30, edge] (d5_10) to (d5_2);
        \draw[bend right=30, edge] (d5_11) to (d5_2);
        \draw[bend right=30, edge] (d5_12) to (d5_2);
        \draw[bend right=30, edge] (d5_2) to (d5_7);
        \draw[bend right=30, edge] (d5_2) to (d5_8);
        \draw[bend right=30, edge] (d5_2) to (d5_9);
        \draw[bend right=30, edge] (d5_7) to (d5_3);
        \draw[bend right=30, edge] (d5_8) to (d5_3);
        \draw[bend right=30, edge] (d5_9) to (d5_3);
        
      \end{tikzpicture}
      \caption{The digraphs $D_{1,0}$, $D_{3,0}$ and $D_{5,0}$.}
      \label{construction_D_k}
    \end{center}
  \end{minipage}
\end{figure}

\begin{proof}
    Let us fix $\ell \in \mathbb{N}$. We define $D_{\ell+1,\ell}$ as the bidirected complete digraph on $\ell+1$ vertices. Note that $D_{\ell+1,\ell}$ clearly satisfies the desired properties.

    Then, for every $k\geq \ell+2$, we iteratively build $D_{k,\ell}$ from $D_{k-1,\ell}$ or $D_{k-2,\ell}$ as follows:
    
\begin{itemize}
    \item If $k+\ell$ is even, we just add a dominating vertex to $D_{k-1,\ell}$ to construct $D_{k,\ell}$. We obtain that $\omega(\UG(D_{k,\ell})) = 1 + \omega(\UG(D_{k-1,\ell})) = k$, $\Delta(B(D_{k,\ell})) = \Delta(B(D_{k-1,\ell})) = \ell$ and $\dic(D_{k,\ell}) = \dic(D_{k-1,\ell}) = \left\lceil \frac{k+\ell-1}{2} \right\rceil = \left\lceil \frac{k+\ell}{2} \right\rceil$ (the last equality holds because $k+\ell$ is even).

    \item If $k+\ell$ is odd (implying that $k$ is at least $\ell+3$), we start from $T$, a copy of $TT_{\frac{k+\ell+1}{2}}$, the transitive tournament on $\frac{k+\ell+1}{2}$ vertices. Note that $\frac{k+\ell+1}{2}\leq k-1$ because $k\geq \ell+3$.
    
    For each arc $xy$ in $T$, we add a copy $D^{xy}$ of $D_{k-2,\ell}$ with all arcs from $y$ to $D^{xy}$ and all arcs from $D^{xy}$ to $x$. Let $D_{k,\ell}$ be the obtained digraph. 
    
    First, $\UG(D_{k,\ell})$ is chordal because it has a perfect elimination ordering: we first eliminate each copy $D^{xy}$ of $D_{k-2,\ell}$, which is possible because $\UG(D_{k-2,\ell})$ is chordal, and $x,y$ are adjacent to every vertex of $D^{xy}$. When every copy of $D_{k-2,\ell}$ is eliminated, the remaining digraph is $T$, which is clearly chordal because it is a tournament.
    
    Next, we have $\omega(\UG(D_{k,\ell})) = \max (\omega(\UG(T)), \omega(\UG(D_{k-2,\ell})) + 2) = k$, and $\Delta(B(D_{k,\ell})) = \Delta(B(D_{k-2,\ell})) = \ell$.
    
    Finally, let us show that $\dic(D_{k,\ell}) \geq \frac{k+\ell+1}{2}$ (the equality then comes from Proposition~\ref{prop:bounded_max_degree}). In order to get a contradiction, assume that $\phi$ is a dicolouring of $D_{k,\ell}$ that uses at most $\frac{k +\ell - 1}{2}$ colours. We know by induction that each copy of $D_{k-2,\ell}$ uses all the colours in $\phi$. Since $T$ is a tournament on $\frac{k+\ell+1}{2}$ vertices, we know that it must contain a monochromatic arc $xy$. Now let $z$ be a vertex in $D^{xy}$ such that $\phi(x)=\phi(y)=\phi(z)$, then $xyz$ is a monochromatic triangle, a contradiction.
\end{itemize}
Figure~\ref{construction_D_k} illustrates the construction of $D_{1,0}$, $D_{3,0}$ and $D_{5,0}$.

\end{proof}

\section{Super-orientations of chordal graphs with a bidirected graph having bounded maximum average degree}
\label{section:bounded_mad}

This section is devoted to the proof of Theorem~\ref{thm:bounded_mad}.
We first need to prove the following.

\begin{lemma}
    \label{lemma-chordal-graph}
    Let $G=(V,E)$ be a chordal graph. There exists an ordering $a_1,\dots ,a_n$ of $V$ such that for any $k\in [n]$ :
    \begin{align*}
         |N(a_k)| &\leq \omega(G)+k-2 ~~~~~~{\rm (P1)}\\
         \text{ and }\ \left | \bigcup_{i=1}^k N[a_i] \right | &\leq \omega(G) + 2k-1  ~~~~~{\rm (P2)}
    \end{align*}
\end{lemma}

\begin{proof}
    Let $(T=(I,F),\mathcal{X} =(B_u)_{u\in I})$ be a valid tree-decomposition of $G$ of width $\omega(G)-1$, which exists by Lemma~\ref{lem:valid_decomposition} (recall that $\tw(G) = \omega(G)-1$ by Proposition~\ref{prop:treewidth_chordal}). One can easily show that, since $T$ is valid, $|I| = n-\omega(G)+1$ (see~\cite[Lemma~2.5]{bodlaenderSIAMJC25}).
    
    Let $P=u_{0},\dots ,u_{r}$ be a longest path in $T$. We root $T$ in $u_r$. For any vertex $u$ of $T$ different from $u_r$, $\father(u)$ denotes the father of $u$ in $T$.
    
    We now consider a Depth-First Search of $T$ from $u_r$. The vertices of $P$ have the priority. Along this route, we label the vertices of $T$. A vertex is labelled when all of its children are labelled. We denote by $v_1,\dots ,v_{n-\omega(G)+1}$ the vertices of $T$ in this labelling. Note that $v_1$ corresponds to $u_0$ and $v_{n-\omega(G)+1}$ corresponds to $u_r$.
    
    Now, for each $i \in \{ 1,\dots , n-\omega(G)\}$, we denote by $a_i$ the unique vertex of $G$ that belongs to $B_{v_i}$ but not to $\father(B_{v_i})$ (recall that $T,\mathcal{X}$ is valid so $a_i$ is well defined). We finally label $a_{n-\omega(D)+1},\dots ,a_n$ the remaining vertices of $G$ in $B_{u_r}$ in an arbitrary way.
    See Figure~\ref{fig-lemme-cordal} for an example of building $a_1,\dots,a_n$.

\begin{figure}[hbtp]
  \begin{minipage}{\linewidth}
    \begin{center}	
      \begin{tikzpicture}[thick,scale=.7, every node/.style={transform shape}]
	  	\tikzset{vertex/.style = {shape=circle,draw, minimum size=2em}}
	  	\tikzset{edge/.style = {->,> = latex'}}
	  
        \node[vertex] (a) at  (-1,3) {a};
        \node[vertex] (b) at  (-1,-1) {b};
        \node[vertex] (c) at  (0,2) {c};
        \node[vertex] (d) at  (0,0) {d};
        \node[vertex] (e) at  (2,2) {e};
        \node[vertex] (f) at  (2,0) {f};
        \node[vertex] (g) at  (3.732,1) {g};
        \node[vertex] (h) at  (3.732,3) {h};
        \node[vertex] (i) at  (2,4) {i};
        \node[vertex] (j) at  (5.732,2) {j};
        \node[vertex] (k) at  (3.732,-1) {k};
        \node[vertex] (l) at  (2,-2) {l};
        \node[vertex] (m) at  (5.732,0) {m};
        
        \draw[] (a) to (d);
        \draw[] (a) to (e);
        \draw[] (a) to (c);
        \draw[] (b) to (d);
        \draw[] (b) to (c);
        \draw[] (b) to (f);
        \draw[] (c) to (d);
        \draw[] (c) to (e);
        \draw[] (c) to (f);
        \draw[] (f) to (e);
        \draw[] (e) to (d);
        \draw[] (d) to (f);
        \draw[] (g) to (d);
        \draw[] (g) to (e);
        \draw[] (g) to (f);
        \draw[] (h) to (g);
        \draw[] (h) to (e);
        \draw[] (h) to (f);
        \draw[] (i) to (e);
        \draw[] (i) to (g);
        \draw[] (i) to (h);
        \draw[] (j) to (e);
        \draw[] (j) to (g);
        \draw[] (j) to (h);
        \draw[] (k) to (g);
        \draw[] (k) to (e);
        \draw[] (k) to (f);
        \draw[] (l) to (f);
        \draw[] (l) to (g);
        \draw[] (l) to (k);
        \draw[] (m) to (f);
        \draw[] (m) to (g);
        \draw[] (m) to (k);
        
        \node[vertex] (acde) at  (8,2.5) {acde};
        \node[] at (8.8,2.5) {$B_1$};
        \node[vertex] (bcdf) at  (8,-0.5) {bcdf};
        \node[] at (8.8,-0.5) {$B_2$};
        \node[vertex] (cdef) at  (9,1) {cdef};
        \node[] at (8.2,1) {$B_3$};
        \node[vertex] (defg) at  (11,1) {defg};
        \node[] at (11.8,1) {$B_7$};
        \node[vertex] (efgh) at  (12,2.5) {efgh};
        \node[] at (11.2,2.5) {$B_6$};
        \node[vertex] (efgk) at  (12,-0.5) {efgk};
        \node[] at (11.2,-0.5) {$B_9$};
        \node[vertex] (eghi) at  (11,4) {eghi};
        \node[] at (11.8,4) {$B_4$};
        \node[vertex] (fgkl) at  (11,-2) {fgkl};
        \node[] at (11.8,-2) {$B_8$};
        \node[vertex] (eghj) at  (14,2.5) {eghj};
        \node[] at (14.8,2.5) {$B_5$};
        \node[vertex] (fgkm) at  (14,-0.5) {fgkm};
        \node[] at (14.9,-0.5) {$B_{10}$};
        
        \draw[edge,orange, dashed] (cdef) to (acde);
        \draw[edge] (cdef) to (bcdf);
        \draw[edge,orange, dashed] (defg) to (cdef);
        \draw[edge] (defg) to (efgh);
        \draw[edge,orange, dashed] (efgk) to (defg);
        \draw[edge] (efgh) to (eghi);
        \draw[edge] (efgh) to (eghj);
        \draw[edge,orange, dashed] (fgkm) to (efgk);
        \draw[edge] (efgk) to (fgkl);
      \end{tikzpicture}
      \caption{A chordal graph $G$ (on the left) and its valid tree-decomposition $T$ (on the right). The orange dashed arcs represent the chosen maximum path $P$. The ordering $a_1,\dots ,a_n$ of $V(G)$ we built is $a,b,c,i,j,h,d,l,e,f,g,k,m$.}
      \label{fig-lemme-cordal}
    \end{center}    
  \end{minipage}
\end{figure}
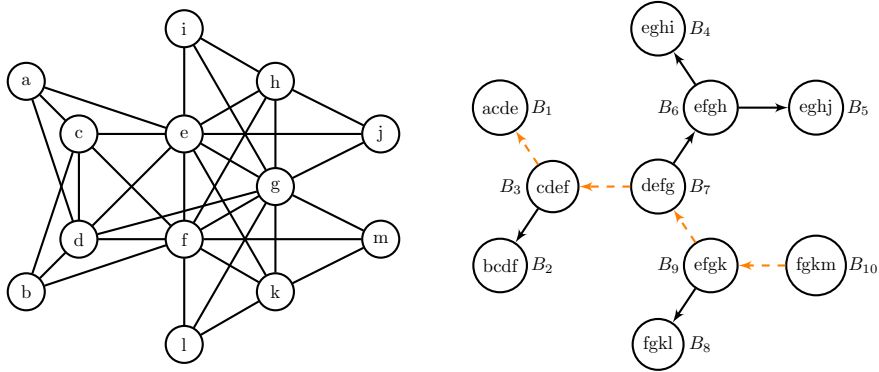
    
    We will now prove that $(a_i)_{1 \leq i \leq n}$ satisfies the two properties of the statement. First observe that, for every $i\in [n]$, $N(a_i) \subseteq \{a_1,\dots,a_{i-1}\} \cup X_{v_i}$ because $a_i\notin \bigcup_{j=i+1}^{n-\omega(G)+1}X_{v_j}$. Hence we have $|N(a_i)|\leq i-1 + \omega - 1 = \omega(G) - 2 + i$, which shows (P1).
    
    To show that (P2) holds, we fix $k\in [n]$. Note that the result is trivially true when $k\geq n-\omega + 1$, thus we assume that $k\leq n-\omega$. Hence, both $v_k$ and $\father(v_k)$ are well defined. We set $X_T = \{v_1,\dots ,v_k\}$, $X_G=\{a_1,...,a_k\}$ and we let $T'$ be the smallest subtree of $T$ that contains all vertices of $X_T$. Let $\ell$ be the largest integer such that $u_{\ell}$ belongs to $V(T')$ ($\ell$ is well defined  because $T'$ contains $v_1=u_0$). We root $T'$ in $u_{\ell}$.
    
    We will now show that $T'$ contains at most $2k$ vertices. If $u_{\ell}=v_k$, then the vertices of $T'$ are exactly $\{v_1,\dots ,v_k\}$ and this is clear. Otherwise let us show that $T''=T'\setminus X_T$ contains at most $k$ vertices, and we will get the result since $|X_T| = k$. By construction we know that every descendant of a vertex $v_i$ is labelled less than $i$. Hence, $T''=T'\setminus X_T$ is a tree rooted in $u_{\ell}$. 
    
    Assume first that $T''$ contains at least two leaves $f_1$ and $f_2$ different from $u_\ell$ ($u_\ell$ may be a leaf it has only one child). We denote by $P_1$ and $P_2$ two paths from their lowest common ancestor. Without loss of generality, we assume that $f_1$ is before $f_2$ in $(v_1,\dots ,v_n)$. Since $f_2$ has a child $g_2$ in $X_T$ and by construction of $(v_i)_{1 \leq i \leq n}$, the internal vertices of $P_1$ are before $g_2$ in $(B_1,\dots ,B_n)$. This implies that all internal vertices in $P_1$ must belong to $X_T$, which contradicts the existence of $f_1$. This shows that $T''$ must have exactly two leaves (one of them is $u_{\ell}$) and then $T''$ is a path rooted in $u_{\ell}$. Since $P$ is a longest path in $T$, we get that $|V(T'')|\le \ell \le k$ and $T'$ contains at most $2k$ vertices as desired.

    \medskip

    We now consider the set $N_G= \{ a_j \in V(G) \mid v_j \in V(T') \setminus \{ u_\ell \} \}$.
    Let $x$ be any vertex in $X_G$. Then every neighbour of $x$ must belong to some bag in $T'$. Moreover, if a vertex belongs to a bag of $T'$, then either it belongs to $B_{u_\ell}$ or it belongs to $N_G$. Then the neighbourhood of $x$ is a subset of $N_G \cup B_{u_\ell}$. Also, $x$ itself belongs to $N_G$. Since $x$ is any vertex in $X_G$, we have: $$\bigcup_{x \in X_G} N[x] \subseteq ( N_G \cup B_{u_\ell} )$$
    
    Since $|N_G| \leq 2k-1$ and $|B_{u_\ell}| = \omega-1$, we get (P2).
\end{proof}

In order to prove Theorem~\ref{thm:bounded_mad}, we prove the more general following result.

\begin{theorem}\label{theo-cordal-madk}
    Let $D$ be a super-orientation of a chordal graph $G$ such that $\Mad(B(D)) \leq d$. For every $\epsilon > 0$, we have
    \[
     \vec{\chi}(D) \leq \left (\frac{1+\epsilon}{2}\right )\omega(G) + \frac{d}{\epsilon} + 1\\
    \]
\end{theorem}
\begin{proof}
    Let  $\epsilon > 0$ and $d\geq 1$, we assume that $\epsilon \leq 1$ for otherwise the result is trivial. We fix $c_{d,\epsilon} = \max\left(\left \lceil \frac{d}{2\epsilon} \right \rceil, \frac{3}{4}d+\frac{d}{8\epsilon} + \frac{1}{2}\right)$. Straightforward calculations imply $c_{d,\epsilon} \leq \frac{d}{\epsilon}+1$. We will show that every super-orientation $D$ of a chordal graph $G$ with $\Mad(B(D))\leq d$ satisfies
    \[ \dic(D) \leq \left(\frac{1+\epsilon}{2}\right)\omega(G) + c_{d,\epsilon}\]

    We prove it by reductio ad absurdum, so assume that $D=(V,A)$ is a smaller counterexample, meaning that $\dic(D) > \left(\frac{1+\epsilon}{2}\right)\omega(G) + c_{d,\epsilon}$. Thus $D$  must be vertex-dicritical (meaning that $\dic(H) < \dic(D)$ for every induced subdigraph $H$ of $D$), for otherwise there exists a vertex $x\in V$ such that $\dic(D-x) = \dic(D)$, and $D-x$ would be a smaller counterexample. 

    For the simplicity of notations, from now on, we write $\omega$ for $\omega(G)$  .
    Let $v$ be any vertex of $D$ and $\alpha$ be any optimal dicolouring of $D-v$ (meaning that $\alpha$ uses exactly $\dic(D) - 1$ colours). Then $\alpha$ cannot be extended to $D$ without using a new colour for $v$ (because $D$ is dicritical). Since every digon (incident to $v$) may forbid at most one colour at $v$, and each pair of simple arcs (incident to $v$) may forbid at most one colour at $v$, we get the following inequalities with $\dig(v)$ the number of digons incident to $v$:
    
    \begin{align}
        \dig(v) + \frac{|N(v)| - \dig(v)}{2} &\geq \dic(D) - 1 >  \left(\frac{1+\epsilon}{2}\right)\omega + c_{d,\epsilon} - 1 \\
        {\rm implying ~~~~~~~~~~~} \dig(v) &> (1+\epsilon)\omega + 2c_{d,\epsilon} - 2 - |N(v)| \label{in-dig}
    \end{align}
    Note that these inequalities hold for every vertex $v$ of $D$. By Lemma~\ref{lemma-chordal-graph}, there is an ordering $a_1, \dots, a_n$ of $V(D)$ such that, for any $i\in [n]$,
    \begin{align*}
         |N(a_i)| &\leq \omega + i -2 ~~~~~~{\rm (P1)}\\
         \text{ and }\ \left | \bigcup_{j=1}^i N(a_j) \right | &\leq \omega + 2i -1 ~~~~~{\rm (P2)}
    \end{align*}
    
    Let us fix $i = \left \lceil \frac{d}{2\epsilon} \right \rceil$. Note that $i \leq c_{d,\epsilon}$. Thus, since $\dic(D) > c_{d,\epsilon}$, we obviously have $i\leq n$. Let $X=\{a_j \mid j\leq i\}$ and $W= \bigcup_{j=1}^i N[a_j]$.
    Together with inequality~\eqref{in-dig}, property (P1) implies, for every $j\in [i]$, $\dig(a_j) > \epsilon\omega + 2c_{d,\epsilon} - j$. Hence we get:
    \begin{equation}
      \sum_{v\in X}\dig(v) = \sum_{j=1}^i\dig(a_j) > \epsilon\omega i + 2c_{d,\epsilon}i - \frac{i(i+1)}{2} \label{eq:lower_bound}
    \end{equation}

    By (P2), we know that $|W|\leq \omega + 2i-1$. Thus $D\ind{W}$ contains at most $\frac{d}{2}(\omega + 2i - 1)$ digons.
    Similarly, since $|X| = i$, $D\ind{X}$ contains at most $\frac{di}{2}$ digons. When we sum $\dig(v)$ over all vertices $v$ in $X$, we count exactly once every digon between $X$ and $W\setminus X$, and exactly twice every digon in $X$. Then, the following is a consequence of~(\ref{eq:lower_bound}).
    
    \begin{eqnarray*}
    \epsilon\omega i + 2c_{d,\epsilon}i - \frac{i(i+1)}{2} < \sum_{v\in X} \dig(v) & \leq & \dig(D\ind{W}) + \dig(D\ind{X}) \\
    & \leq & \frac{d}{2}(\omega + 2i - 1) + \frac{di}{2}\\
    \end{eqnarray*}

    Since $i=\left \lceil \frac{d}{2\epsilon} \right \rceil$, we conclude that $c_{d,\epsilon} < \frac{3}{4}d+\frac{d}{8\epsilon} + \frac{1}{2}$, a contradiction.
\end{proof}

The proof of Theorem~\ref{thm:bounded_mad} now follows.

\boundedmad*

\begin{proof}
    This is a direct consequence of Theorem~\ref{theo-cordal-madk} applied for $\epsilon = \sqrt{\frac{d}{\omega(G)}}$.
\end{proof}

\section{Super-orientations of chordal graphs \texorpdfstring{$\bid{C_4}$}{C4}-free}
\label{section:no_C4}

This section is devoted to the proof of Theorems~\ref{thm:C4_free} and~\ref{construction:C4_free}.

\Cfourfree*

\begin{proof}
    We assume that $\omega = \omega(G)$ is odd, otherwise we select an independent set $I$ of $D$ such that $D' = D-I$ satisfies $\omega(\UG(D')) = \omega - 1$, so $\omega(\UG(D'))$ is odd and $\dic(D) \leq \dic(D') + 1$ (the existence of $I$ is guaranteed because $G$ is chordal).

    Let $(T,\mathcal{X} =(B_u)_{u\in V(T)})$ be a valid tree-decomposition of $G$, that is each bag $B\in \mathcal{X}$ has size exactly $\omega$ and, for every two adjacent bags $B$ and $B'$, $|B \setminus B'| = 1$. Recall that the existence of such a tree-decomposition is guaranteed by Lemma~\ref{lem:valid_decomposition}. We assume that each bag induces a clique on $G$, otherwise we just add the missing arcs (oriented in an arbitrary direction). Note that this operation does increase $\omega$ nor decrease $\dic(D)$ and does not create any $\bid{C_4}$.
    
    Let $k=\frac{\omega+3}{2}$. A $k$-dicolouring $\phi$ of $D$ is \textit{balanced} if, for each bag $B$ and colour $c\in [k]$, $0\leq |\phi^{-1}(c) \cap B| \leq 2$. Note that every balanced $k$-dicolouring satisfies $|\phi^{-1}(c) \cap B|=1$ for either 1 or 3 colours. Moreover, in the former case, exactly one colour of $[k]$ is missing  in $\phi(B)$.
    We will show that $\dic(D) \leq k$ by proving the existence of a balanced $k$-dicolouring $\phi$ of $D$ such that, for each bag $B$, we have:
    \begin{itemize}
        \item[(1)] $|\phi^{-1}(c) \cap B|=1$ holds for exactly one colour $c$, or
        
        \item[(2)] $|\phi^{-1}(c_i) \cap B|=1$ for exactly three distinct colours $c_1,c_2,c_3$ and two vertices of $\{v_1, v_2, v_3\}$ are connected by a $\bid{P_3}$ in $D$ (where $\{v_i\} = \phi^{-1}(c_i)\cap B$ and a $\bid{P_3}$ is a bidirected path on 3 vertices). 
    \end{itemize}
    We will say that a bag $B$ is of type (1) or (2), depending if $\phi$ satisfies condition (1) or (2) respectively on $B$.

    We show the existence of $\phi$ by induction on the number of bags in the tree-decomposition. If $|V(T)| = 1$, let $\mathcal{X} = \{B\}$, then $D$ is a semi-complete digraph on $\omega$ vertices which is $\bid{C_4}$-free. We construct $\phi$ greedily as follows: choose a simple arc $uv$ such that both $u$ and $v$ have not been coloured yet, and use a new colour for them. At the end, there are either one or three uncoloured vertices. If there is only one, we just use a new colour for it and $B$ is of type (1), otherwise the three remaining vertices induce a bidirected triangle on $D$ and we can use one new colour for each of them, so $B$ is of type (2).

    Assume now that $|V(T)| \geq 2$. Let $x$ be a leaf of $T$ and $y$ its only neighbour in $T$. Let $\{u\}= B_y \setminus B_x$ and $\{v\}= B_x \setminus B_y$. By induction, with $D-v$ and $(T-x, \mathcal{X}\setminus B_x)$ playing the role of $D$ and $(T,\mathcal{X})$ respectively, there exists a balanced $k$-dicolouring $\phi$ of $D-v$ for which each bag is of type (1) or (2). We will show by a case analysis that $\phi$ can be extended to $v$.
    \begin{itemize}
        \item Assume first that $B_y$ is of type (1), and let $r$ be the only vertex alone in its colour class in $D\ind{B_y}$. If $r=u$, then we set $\phi(v) = \phi(u)$ and $\phi$ is a balanced $k$-dicolouring of $D$ with $B_x$ being of type (1). Henceforth assume $u\neq r$. Let $w$ be the neighbour of $u$ in $B_y$ such that $\phi(w) = \phi(u)$. Since $u$ and $v$ are not adjacent, setting $\phi(v) = \phi(u)$ yields a balanced $k$-dicolouring of $D$, with $B_x$ being of type (1),  except if $w$ and $v$ are linked by a digon. Analogously, setting $\phi(v) = \phi(r)$ yields a balanced $k$-dicolouring of $D$, with $B_x$ being of type (1) since $|\phi^{-1}(c) \cap B_x|=1$ holds only for $c=\phi(w)$,  except if $r$ and $v$ are linked by a digon.

        But then, if both $[v,w]$ and $[v,r]$ are digons, we can set $\phi(v)$ to the missing colour of $\phi(B_y)$. Then $\phi$ is a balanced $k$-dicolouring of $D$ with $B_x$ being of type (2), since $|\phi^{-1}(c) \cap B_x|=1$ holds exactly for every $c\in \{\phi(w), \phi(v),\phi(r) \}$ with $r,w$ being connected by a $\bid{P_3}$ in $D$.

        \item Henceforth assume that $B_y$ is of type (2) and let $r,s,t$ be the only vertices alone in their colour class in $D\ind{B_y}$ such that $s$ and $t$ are connected by a $\bid{P_3}$ in $D-v$. If $u=r$, then we set $\phi(v) = \phi(u)$ and $\phi$ is a balanced $k$-dicolouring of $D$ with $B_x$ being of type (2).

        Assume now that $u\in \{s,t\}$. Without loss of generality, we assume that $u=s$. If $r$ and $v$ are not linked by a digon, we can set $\phi(v) = \phi(r)$ and $\phi$ is a balanced $k$-dicolouring of $D$ with $B_x$ being of type (1). The same argument holds if $t$ and $v$ are not linked by a digon. But if both $[v,r]$ and $[v,t]$ are digons, we can set $\phi(v) = \phi(s)$. Then $\phi$ is a balanced $k$-dicolouring of $D$ with $B_x$ being of type (2), since $|\phi^{-1}(c) \cap B_x|=1$ holds exactly for every $c\in \{\phi(v), \phi(r),\phi(t) \}$ with $r,t$ being connected by a $\bid{P_3}$ in $D$.

        Assume finally that $u\notin \{r,s,t\}$ and let $w$ be the neighbour of $u$ in $B_y$ such that $\phi(w) = \phi(u)$. If $r$ and $v$ are not linked by a digon, we can set $\phi(v) = \phi(r)$ and $\phi$ is a balanced $k$-dicolouring of $D$ with $B_x$ being of type (2), where $|\phi^{-1}(c) \cap B_x|=1$ holds exactly for every $c\in \{\phi(w), \phi(s),\phi(t) \}$ with $s,t$ being connected by a $\bid{P_3}$ in $D-v$. The same argument holds if $v$ and $w$ are not linked by a digon. Henceforth we assume that both $[v,w]$ and $[v,r]$ are digons. Since $D$ is $\bid{C_4}$-free, and because $s,t$ are connected by a $\bid{P_3}$ in $D-v$, we know that either $[v,s]$ or $[v,t]$ is not a digon of $D$. Assume without loss of generality that $[v,s]$ is not, then we set $\phi(v) = \phi(s)$.  Then $\phi$ is a balanced $k$-dicolouring of $D$ with $B_x$ being of type (2), since $|\phi^{-1}(c) \cap B_x|=1$ holds exactly for every $c\in \{\phi(w), \phi(r),\phi(t) \}$ with $w,r$ being connected by a $\bid{P_3}$ in $D$.
    \end{itemize}
\end{proof}

\constructionCfourfree*

\begin{proof}
    We only have to prove it for $k=3$. For larger values of $k$, we build $D_{k,n}$ from $D_{k-1,n}$ or $D_{k-2,n}$ as in the proof of Theorem~\ref{construction_chordal}. 
    The digraph $D_{3,n}$, depicted in Figure~\ref{fig:C4_free_bound}, is clearly a super-orientation of a $2$-tree. As a consequence of Theorem~\ref{thm:hajos_join}, it has dichromatic number 3, since it is obtained from successive Haj\'os joins applied on $\bid{K_3}$.
\end{proof}
    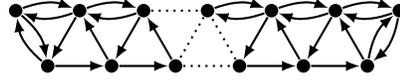
\begin{figure}[hbtp]
      \begin{minipage}{\linewidth}
        \begin{center}	
          \begin{tikzpicture}[thick,scale=.85, every node/.style={transform shape}]
            \tikzset{vertex/.style = {circle,fill=black,minimum size=6pt, inner sep=0pt}}
            \tikzset{edge/.style = {->,> = latex}}
                
            \foreach \i in {0,...,5}{
                \node[vertex] (u\i) at  (\i,0) {};
                \node[vertex,shift={(\i,0)}] (v\i) at  (-60:1) {};
            }
            \node[vertex] (u6) at  (6,0) {};
            \foreach \i in {0,...,1}{
                \pgfmathtruncatemacro{\j}{\i + 1}
                \draw[edge,bend left=20] (u\j) to (u\i){};
                \draw[edge,bend left=20] (u\i) to (u\j){};
                \draw[edge] (v\i) to (v\j){};
            }
            
            \draw[edge,bend left=20] (u0) to (v0){};
            \draw[edge,bend left=20] (v0) to (u0){};
            \draw[edge] (u1) to (v0){};
            \draw[edge] (u2) to (v1){};
            \draw[edge] (v1) to (u1){};
            \draw[edge] (v2) to (u2){};
            
            \draw[dotted] (u2) -- (u3) -- (v2) -- (v3) --(u3);
            
            \foreach \i in {3,...,4}{
                \pgfmathtruncatemacro{\j}{\i + 1}
                \draw[edge,bend left=20] (u\j) to (u\i){};
                \draw[edge,bend left=20] (u\i) to (u\j){};
                \draw[edge] (v\i) to (v\j){};
            }
            \draw[edge] (u4) to (v3){};
            \draw[edge] (u5) to (v4){};
            \draw[edge] (v4) to (u4){};
            \draw[edge] (v5) to (u5){};
            \draw[edge,bend left=20] (u6) to (v5){};
            \draw[edge,bend left=20] (v5) to (u6){};
            \draw[edge,bend left=20] (u5) to (u6){};
            \draw[edge,bend left=20] (u6) to (u5){};
          \end{tikzpicture}
          \caption{The digraph $D_{3,n}$.}
          \label{fig:C4_free_bound}
        \end{center}
      \end{minipage}
    \end{figure}

\section{Further research}

In this work, we gave both lower and upper bounds on the dichromatic number orientations and super-orientations of different classes of chordal graphs and cographs. A lot of open problems arise and we detail a few of them.
First, we do not know if the bound of Theorem~\ref{thm:bounded_mad} is optimal, and we ask the following.

\begin{problem}
    Show the existence of $f$ such that every super-orientation $D$ of a chordal graph $G$ satisfies $\dic(D) \leq \frac{1}{2}\omega(G) + f(\Mad(B(D)))$.
\end{problem}

We also ask if Theorem~\ref{thm:C4_free} is true not only for $\bid{C_4}$-free digraphs but for $\bid{C_\ell}$-free digraphs in general.

\begin{problem}
Show the existence of $g$ that every $\bid{C_\ell}$-free super-orientation $D$ of a chordal graph $G$ satisfies $\dic(D) \leq \frac{1}{2}\omega(G) + g(\ell)$.
\end{problem}

A famous class of graphs is the class of claw-free graphs (a graph is {\it claw-free} if it does not contain $K_{1,3}$ as an induced subgraph). Line-graphs and proper interval graphs are examples of claw-free graphs. We ask the following.

\begin{problem}
    Does every orientation $\Vec{G}$ of a claw-free graph $G$ satisfy $\dic(\Vec{G}) = O\left(\frac{\omega(G)}{\log \omega(G)}\right)$?
\end{problem}

A celebrated conjecture of Erd\H{o}s and Neumann-Lara (see~\cite{erdosPNCN1979}) states that every orientation $\Vec{G}$ of a graph $G$ satisfies $\dic(\Vec{G}) = O\left(\frac{\Delta(G)}{\log \Delta(G)}\right)$. Since every claw-free graph $G$ satisfies $\Delta(G) \leq 2\omega(G)-2$, Erd\H{o}s and Neumann-Lara's conjecture, if true, would answer the question above positively.

\section*{Acknowledgments}

The authors are thankful to Virginia Ardévol Martínez for late discussions on Theorem~\ref{thm:construction_interval}.

\end{document}